\documentstyle[12pt]{article}
\begin{document}
\begin{center}
Applications of Generalized Special Functions in Stellar
Astrophysics\\
\vspace{0.5cm}
H.J. Haubold and A.M. Mathai\\
\vspace{0.5cm}
Department of Mathematics and Statistics\\
McGill University\\
Montreal, Canada H3A 2K6\\
\end{center}
\vspace{3cm}
{\bf ABSTRACT}\par
\vspace{0.7cm}
This article gives an brief outline of the applications of
generalized
special functions such as generalized hypergeometric functions, 
G-functions and H-functions into the general area of nuclear energy
generation and reaction rate theory such as the energy generation
in a simple stellar model and nuclear reaction rates in 
non-resonant and resonant as well as screened non-resonant
situations.
\section{Introduction}
During the last few years several authors have been trying to get
explicit and exact analytic representations for various problems in
astropyhysics. Analytical solutions are often extremely difficult
to obtain and often lead to complicated intractable mathematical
expressions. Due to these facts some researchers  go for
approximations, numerical results and computer printouts. But such
numerical results usually do not give much insight into what is
happening between the stage of the formulation of the problem and
the final computer printout. When the experimental results do not
agree with the computer printouts resulting from a proposed theory
one needs explicit analytic formulation to modify the theory. Hence
there has been renewed interest in getting explicit and exact
theoretical results rather than numerical results. During the past
few years the authors have demonstrated that the
techniques of generalized special functions combined with
statistical techniques are very powerful in getting analytical
results for many problems which were believed to be
intractable mathematically. This paper will confine itself to the
application of such techniques for studying energy generation and
thermonuclear reaction rates for stellar models.
\section{Energy Generation Rate in Simple Stellar Models}
Some simple stellar models for looking into the internal structure
of stars have been given by Kourganoff
(1980) and Haubold and Mathai (1984, 1986b). Consider a spherically
symmetric star in quasi-static equilibrium, that is, we exclude
magnetic fields, rotation and other effects which are changing
rapidly, so that we can take the star as a gas sphere in
hydrostatic
equilibrium. Let $r$ denote the distance from the center to an
interior point and $R$ the radius of the star under consideration.

\subsection{Linear model}
The simplest model one can think of for the radial density $\rho
(r)$ is the following,
\begin{equation}
\rho(r)=\rho_c(1-\frac{r}{R})=\rho_c(1-x),
x=r/R,\rho(0)=\rho_c,\rho(R)=0.
\end{equation}
In this case the mass
\begin{equation}
M(r)=4\pi\int_0^rt^2\rho(t)dt=\frac{4\pi}{3}\rho_cr^3(1-
\frac{3x}{4}),
\end{equation}
the pressure
\begin{eqnarray}
P(r)& = & P(0)-G\int_0^r\frac{M(t)}{t^2}\rho(t)dt \nonumber \\
& = &\frac{\pi}{36}G\rho^2_cR^2(5-24x^2+28x^3-9x^4),
\end{eqnarray}
where $G$ is the gravitational constant, and the temperature
from the equation of state $P(r)=k\rho (r)T(r)/M_u\mu,$ where $\mu$
is the mean
molecular weight, $M_u$ is the atomar mass unit, $k$ is the
Boltzmann constant, is then
\begin{eqnarray}
T(r)& = &\frac{\pi G\mu M_u}{36k}\rho_cR^2\frac{(5-24x^2+28x^3-
9x^4)}{1-x}\nonumber \\
& = & \frac{\pi G\mu M_u}{36k}\rho_cR^2\sum_{n=0}^\infty = c_nx^n
\end{eqnarray}
with $c_0=c_1,c_2=-19/5, c_3=9/5, c_3=9/5, c_n=0$ for $n\geq 4.$
The
nuclear energy generation rate $\epsilon (r)$ will depend on $T(r)$
and
$\rho(r)$. If we take a simple model,
\[\epsilon(r)=\epsilon_0(\rho_0, T_0)(\rho(r)/\rho_0)^\alpha
(T(r)/T_0)^\beta, \alpha > 0, \beta > 0\]
then
\begin{equation}
\epsilon(r)= d_1(1-x)^\alpha(\sum^\infty _{n=0}c_nx^n)^\beta,
\end{equation}
where $d_1$ is a constant. The luminosity $L$ of the star is then,
\begin{eqnarray}
L & = & 4\pi\int_0^R r^2\rho(r)\epsilon(r)dr\nonumber \\
& = & d_2\int_0^1 x^2(1-x)^{\alpha +1}(\sum^\infty _{n=0}
c_nx^n)^\beta dx,
\end{eqnarray}
where $d_2$ is a constant. We consider a general integral
\begin{equation}
g(\alpha, \beta, \gamma)=\int_0^1 x^\gamma (1-
x)^{\alpha+1}(\sum^\infty _{n=0} c_nx^n)^\beta dx.
\end{equation}
The $L$ in (6) is evaluated in Haubold and Mathai (1984) for
various
structures of $\sum_{n=0}^\infty = c_nx^n$. One such structure is
when
\[\sum_{n=0}^\infty c_nx^n=(1+a_1x)(1+a_2x)\ldots(1+a_kx)\]
for a fixed $k$. In this case one can show that $g(\alpha, \beta,
\gamma)$ can be evaluated in terms of a Lauricella function $F_D$.
Then a particular case will go in terms of Gauss' hypergeometric
function.
\subsection{A general stellar model}
A linear decrease of the density from the center to the surface may
not be very appropriate. Faster or slower decrease can be
incorporated by using one more parameter and taking the model for
the density as,
\begin{eqnarray}
\rho (r) & = & \rho_c(1-(\frac{r}{R})^\delta)\nonumber \\ 
& = & \rho_c(1-x^\delta), x=r/R, \rho (0)=\rho_c, \rho (R)=0.
\end{eqnarray}
In this case the mass, pressure, and
luminosity will work out to be the following.
\begin{equation}
M(r)=\frac{4\pi}{3}=\rho_cR^3x^3(1-\frac{3}{\delta+3}x^\delta).
\end{equation}
\begin{equation}
P(r)=\frac{4\pi G}{3}\rho_c^2 R^2 \left\{\psi-
\frac{x^2}{2}+\frac{(\delta+6)}{(\delta+2)(\delta+3)}x^{\delta+2}-
\frac{3}{2(\delta+1)(\delta+3)}x^{2\delta+2}\right\},
\end{equation}
\[\psi=\frac{1}{2}-
\frac{\delta+6}{(\delta+2)(\delta+3)}+\frac{3}{2(\delta+1)(\delta
+3)},\]
and
\begin{eqnarray}
L(R)& = & 
4\pi\rho_c\epsilon_0(\rho_c,T_c)\frac{R^3}{\delta}\int_0^1x^{(3/
\delta)-1}(1-x)^{\alpha+\beta+1}\nonumber \\
& & [1-\frac{x^{2/\delta}}{2\psi}\left\{1-
\frac{2(\delta+6)x}{(\delta+2)(\delta+3)}+\frac{3x^2}{(\delta+1)(
\delta+3)}\right\}]^\beta dx.
\end{eqnarray}
For $\beta=1$ or a positive integer the evaluation of the integral
in (11) is not difficult. For the general case one can look at
various situations as in Section 2.1 and solve many of these in
terms of generalized functions.
\subsection{A further generalized stellar model}
A model with more flexibility for the density is a two parameter
model
\[\rho(r)=\rho_c[1-(\frac{r}{R})^\delta]^\gamma, \delta > 0,
\gamma > 0.\]
In this case the mass will lead to the form
\begin{eqnarray}
M(r) & = & 4\pi \rho_c \int_0^r
t^2\left[1-(\frac{t}{R})^\delta \right]^\gamma
dt, 0 < \frac{t}{R}\leq \frac{r}{R}\leq 1 \nonumber \\
& = & \frac{4}{3}\pi\rho_cr^3\; _2F_1\left(-\gamma,
\frac{3}{\delta};\frac{3}{\delta}+1;(\frac{r}{R})^\delta \right),
\end{eqnarray}
where $_2F_1$ is a Gauss' hypergeometric function. Pressure and
temperature in this case can be obtained in series of
hypergeometric functions. For example,
\begin{eqnarray}
P(r) & = &
\frac{4\pi G\rho^2_c}{\delta^2}\left\{R^2\sum_{m=0}^\infty \frac{(-
\gamma)_m(\frac{\delta}{2}+1)_m}{m!(\frac{\delta}{2}+1+\gamma)_m(
\frac{3}{\delta}+m)(\frac{2}{\delta}+m)}\right.\nonumber \\
& - & r^2\sum_{m=0}^\infty \frac{(-
\gamma)_m\left[(\frac{r}{R})^\delta\right]^m}{m!(\frac{3}{\delta}
+m)(\frac{2}
{\delta}+m)} \\
& _2F_1 & \left. \left(-\gamma, \frac{2}{\delta}+m;
\frac{2}{\delta}+m+1;(\frac{r}{R})^\delta\right)\right\}.
\end{eqnarray}
Note that when $\gamma$ is a positive integer then $P(r)$ is
available as a finite sum.
\section{Collision Probabilities}
In Haubold and Mathai (1984a) a detailed derivation of the
collision probability for nuclear reactions is given under the
assumption that the
nuclear velocity distribution remains Maxwell-Boltzmannian. The
relation between the macroscopic nuclear reaction probability and
the microscopic nuclear reaction cross section on the basis of
equilibrium-thermodynamic arguments is studied in this paper as
well as in Haubold and Mathai (1985). A basic integral to be
evaluated to compute the reaction probability, under this approach,
is the following:
\begin{equation}
A=\int_0^\infty y^\nu e^{-y-zy^{-1/2}}dy
\end{equation}
More complicated integrals will appear when dealing with
thermonuclear reaction rates in various non-resonant and resonant
cases, see for example, Haubold and Mathai (1986, 1986a,b,c). First
we will introduce a statistical technique of tackling integrals of
the type (15). Let us consider a slightly more general integral
\begin{equation}
A_r(z;p,n,m)=p\int_0^\infty e^{-py}y^{-nr}e^{-zy^{-n/m}}dy.
\end{equation}
{\bf Theorem 3.1.}\par
\begin{eqnarray}
& p\int_0^\infty e^{-pt}t^{-nr}e^{-zt^{-n/m}}dt & \nonumber \\
& & =p^{nr}(2\pi)^{(2-n-m)/2}m^{1/2}n^{(1-2nr)/2}\\
& & G^{m+n,0}_{0,m+n}(\frac{z^m p^n}{m^m n^n}\mid
_{0,\frac{1}{m},\ldots,\frac{m-1}{m},
\frac{1-nr}{n},\ldots,\frac{n-
nr}{n}}),\nonumber 
\end{eqnarray}
for $R(p) > 0,R(z) > 0,G(.)$ is a G-function and R(.)
denotes the real part of (.). For a discussion of the G-function
and the more generalized H-function see Mathai and Saxena (1973,
1978).\par
Writing the Mellin-Barnes integral representation for the G-
function and then simplifying the gammas by using the
multiplication formula for gamma functions, namely,
\begin{equation}
\Gamma (mz)=(2\pi)^{(1-m)/2} m^{mz-
1/2}\Gamma(z)\Gamma(z+\frac{1}{m})\ldots \Gamma(z+\frac{m-1}{m}),
m=1,2\ldots
\end{equation}
one has the right side of (17) simplified to the following form,
\begin{equation}
A_r(z;p,n,m)=p^{nr}\frac{m}{n}\frac{1}{2\pi i}\int_L
\Gamma(\frac{ms}{n})\Gamma(1-nr+s)(pz^{m/n})^{-s}ds,
\end{equation}
where $L$ is a suitable contour. The integral on the right side of
(19) can be written as an H-function, see Mathai and Saxena
(1978), which can be reduced to a G-function due to the fact that
$m/n$ is rational. Now consider two independent positive real
random variables $x$ and $y$ with the density functions $f_1(x)<
0$ for $0 < x < \infty, f_2(y) < 0$ for $0 < y < 
\infty$ and $f_1(x)=0,f_2(y)=0$ elsewhere. Let the (s-1)st moments
of $x$ and $y$ be denoted by $g_1(s)$ and $g_2(s)$ respectively.
Consider $u=xy$. From statistical independence of $x$ and $y$ one
has
\[E(u^{s-1})=(Ex^{s-1})(Ey^{s-1})=g_1(s)g_2(s),\]
where $E$ denotes the expected value. If the density of $u$ is
denoted by $g(u)$ then from the inverse Mellin transform
\begin{equation}
g(u)=\frac{1}{2\pi i}\int_{l_1}g_1(s)g_2(s)u^{-s}ds,
\end{equation}
where $L_1$ is a suitable contour. Now by using transformation of
variables technique the density is given by
\begin{equation}
g(u)=\int_0^\infty f_1(v)f_2(\frac{u}{v})v^{-1}dv.
\end{equation}
But due to the uniqueness of the density of $u$ the $g(u)$
appearing in (20) and (21) must be one and the same. Let $v=pt,
u=pz^{m/n},f_1(t)=t^{1-nr}e^{-t} \mbox{and} f_2(t)=e^{-t^{n/m}}.$
Then one
has
$$g_1(s)=\Gamma(1-nr+s),R(1-nr+s)> 0$$
$$g_2(s)=\frac{m}{n}\Gamma(\frac{ms}{n}),R(s)> 0$$
and
$$\frac{1}{2\pi i}\int_{L_1}g_1(s)g_2(s)u^{-s}ds=\frac{1}{2\pi
 i}\int_{L_1}(\frac{m}{n})\Gamma(\frac{ms}{n})\Gamma(1-
nr+s)(pz^{m/n})^{-s}ds$$
\begin{equation}
\int_0^\infty f_1(v)f_2(\frac{u}{v})v^{-1}dv=\int_0^\infty e^{-
pt}(pt)^{1-nr}e^{-zt^{n/m}}t^{-1}dt,
\end{equation}
which gives
\begin{equation}
A_r(z;p,n,m)=p^{nr}(\frac{m}{n})\frac{1}{2\pi i}
\int_{L_1}\Gamma(\frac{ms}{n})\Gamma(1-nr+s)(pz^{m/n})^{-s}ds.
\end{equation}
Comparing (23) and (19) the result is established.\par
Thus in order to evaluate the basic collison probability integral
one has to represent the G-function in Theorem 3.1 in computable
forms. Explicit computable representations of the G-function for
various parameter values are given in Haubold and Mathai (1984a).
Computable representation of a general G-function is available in
Mathai and Saxena (1973).
\subsection{Non-resonant reaction rate}
In this case the closed-form representation of the screened nuclear
reaction rate can be evaluated by evaluating the following
integral.
\begin{equation}
I(z;t,a,v,n,m)=\int^\infty_0e^{-ay}y^ve^{-z(y+t)^{-n/m}}dy,
\end{equation}
for $z>0, t>0, a>0, m,n$ positive integers, see Haubold
and Mathai (1986c) for details.\par
The integral in (24) can be evaluated by using the following
lemmas, of these, Lemma 3.1a can be established by going through
the same process as in the proof of Theorem 3.1, but Lemma 3.1b
needs some modifications. In this case one of the random variables
will be having a non-zero density function in the interval [0,d]
and the density function will be zero outside this interval. Apart
from this modification the derivation will remain more or less
parallel. Hence we will list the lemmas and the theorem to follow
without proofs.\par
\noindent
{\bf Lemma 3.1a.} For $a> 0, z> 0, n,m$ positive integers
\begin{eqnarray}
N_1(z;a,v,n,m) & = & \int_0^\infty e^{-ay}y^ve^{-zy^{-n/m}}dy
\nonumber \\
& = & a^{-(v+1)}(2\pi)^{(2-m-n)/2}m^{1/2}n^{v+1/2} \label{line2} \\
& &
G^{m+n,0}_{0,m+n}(\frac{z^ma^n}{m^mn^n}\mid_{0,\frac{1}{m},\ldots
,\frac{m-1}{m},1+\frac{v}{n},\ldots,n+\frac{v}{n}}). \nonumber 
\end{eqnarray}
\noindent
{\bf Lemma 3.1b.} For $z>0, d>0, a>0, m,n$ positive
integers and denoting $N_2(z;d,a,v,n,m)$ by $N_2$,
\begin{eqnarray}
N_2 & = & \int_0^dy^ve^{-ay}e^{-zy^{-n/m}}dy \nonumber \\
& = & (2\pi)^{\frac{1-m}{2}}m^{1/2}n^{-
1}d^{v+1}\sum^\infty_{r=0}\frac{(-
ad)^r}{r!}G^{m+n,0}_{n,m+n}(\frac{z^m}{d^nm^m}
\mid^{(a)}_{(b)}),\label{line2}\\
(b) & = & (\frac{v}{n}+\frac{(r+1)}{n}+\frac{j-1}{n},
j=1,\ldots,n,\frac{j-1}{m}, j=1,\ldots, m), \nonumber \\
(a) & = & (\frac{v}{n}+\frac{r+2}{n}+\frac{j-1}{m},
j=1,\ldots,n). \nonumber
\end{eqnarray}
\noindent
{\bf Theorem 3.2.} For $a>0, z> 0, t> 0, m,n$ positive
integers and $v$ a non-negative integer
\begin{eqnarray*}
\int_0^\infty y^ve^{-ay}e^{-z(y+t)^{-
n/m}}dy & = & t^{v+1}e^{a_1}\sum_{r=0}^v
(^v_r)(-1)^r[N_1(z_1;a_1,v-
r,n,m) \\ 
& & -N_2(z_1;1,a_1,v-r,n,m)],
\end{eqnarray*} 
where $a_1=at, z_1=zt, (^v_r)=v!/[r!(v-r)!],0!=1,N_1(.)$ and
$N_2(.)$
are given in (25) and (26) respectively.\par
\noindent
\subsection{Resonant reaction rates}
Analytic representations of the reaction probability for the
resonant reaction rates are discussed in detail in Haubold and
Mathai (1986). The integral to be evaluated has a more complicated
structure compared to the integrals in Section 3.1. In this case
the
integral to be evaluated is $N_3$ where,
\begin{equation}
N_3=\int_0^\infty \frac{t^v exp[-at-qt^{-n/m}]}{(b-t)^2+g^2}dt,
\end{equation}
for $a>0, q>0, m,n$ positive integers and $v$ a non-negative
integer. This will be evaluated with the help of Lemmas 3.1a, 3.1b,
Theorems 3.1, 3.2 and the following result which will be stated as
a lemma.\par
\noindent
{\bf Lemma 3.2a.} For $(b-t)^2+g^2> 0,$ \par
\begin{equation}
[(b-t)^2+g^2]^{-1}=\int_0^\infty exp\left\{ -[(b-t)^2+g^2]x
\right\} dx.
\end{equation}
Now replacing the denominator on the right side of (27) by the
integral in (28) and evaluating the double integral with the help
of Lemmas 3.1a, 3.1b and Theorems 3.1 and 3.2 we have the following
result.\par
\medskip
\noindent
{\bf Theorem 3.3.}\par
\begin{eqnarray*}
N_3 & = & \sum_{k=0}^\infty \frac{(-
1)^k}{g^2(g^2)^k}\sum_{k_1=0}^{2k}\left(^{2k}_{k_1}\right)(-
1)^{k_1}\\
& &
b^{2k-k_1}a^{-(v+k_1+1)}(2\pi)^{(2-n-m)/2}m^{1/2}n^{v+k_1+1/2}\\
& &
G^{m+n,0}_{0,m+n}\left(\frac{q^ma^n}{m^mn^n}\mid_{0,\frac{1}{m},
\ldots,\frac{m-1}{m},
\frac{(1+v+k_1)}{n},\ldots,\frac{(n+v+k_1)}{n}}\right),
\end{eqnarray*}
where $N_3(.)$ is defined in (27). One can look at structures
which are mathematically more complicated than the one in (27).
For example the denominator in (27) could be replaced by $[(b-
t)^2+g^2]^\gamma$ for $\gamma \geq 1.$ Still the techniques
described above will work. One can replace $t^{-n/m}$ by $(1+t)^{-
n/m}$ and one can also look at a finite range integral going from
0 to some number $d<\infty$ instead of 0 to $\infty$. Integrals
with such modifications can be tackled by using the combined
statistical and generalized special function techniques described
in Sections 2 and 3 of this paper. 
\clearpage
\noindent
\begin{center}
{\bf References}
\end{center}
Korganoff, V. (1980). {\it Introduction to Advanced Astrophysics.}
Reidel Publication, Boston.\par
\medskip
\noindent
Haubold, H.J. and Mathai, A.M. (1984). On the nuclear energy
generation rate in a simple analytic stellar mode. {\it Annalen der
Physik,} {\bf 41(6)}, 372-379.\par
\medskip
\noindent
Haubold, H.J. and Mathai, A.M. (1984a). On nuclear reaction rate
theory. {\it Annalen der Physik}, {\bf 41(6)}, 380-396.\par
\medskip
\noindent
Haubold, H.J. and Mathai, A.M. (1985). The Maxwell-Boltzmannian
approach to the nuclear reaction  rate theory. {\it Fortschr.
Phys.,} {\bf 33}, 623-644.\par
\medskip
\noindent
Haubold, H.J. and Mathai, A.M. (1986). The resonant thermonuclear
reaction rate. {\it J.Math. Phys.}, {\bf 27}, 2203-2207.\par
\medskip
\noindent
Haubold, H.J. and Mathai, A.M. (1986b). Analytic solution to the
problem of nuclear energy generation rate in a simple stellar
model. {\it Astron. Nachr.,} {\bf 307,}9-12.\par
\medskip
\noindent
Haubold, H.J. and Mathai, A.M. (1986c). Analytic results for
screened non-resonant nuclear reaction rates. {\it Astrophysics and
Space Science,} {\bf 127}, 45-53.\par
\medskip
\noindent
Mathai, A.M. and Saxena, R.K. (1973). {\it Generalized
Hypergeometric Functions with Applications in Statistics and
Physical Sciences.} Lecture Notes in Mathematics No. 348, Springer-
Verlag, Heidelberg.\par
\medskip
\noindent
Mathai, A.M. and Saxena, R.K. (1978). {\it The H-function with
Applications in Statistics and Other Disciplines,} Wiley Halsted,
New York.
\end{document}